\input amstex
\documentstyle{amsppt}
\magnification=1200
\pageheight{9.5 true in}

\def\N{\Bbb N}
\def\Pm{P^*_\mu}
\def\Pl{P^*_\lambda}
\def\l{\lambda}

\def\ga{\gamma}
\def\xm{\prod x_i^{-1}}
\def\Vm{V(x)^{-1}}
\def\L{\Lambda}
\def\Ln{\L(n)}
\def\Lt{\L^t}
\def\Ltn{\Lt(n)}
\def\Lot{\L^{1/t}}
\def\Lotn{\Lot(n)}
\def\Z{\Bbb Z}
\def\D{\frak{D}}
\def\te{\tilde\eta}
\def\etm{\eta_{\text{max}}}
\def\Ts{{\sssize T}(s)}
\def\fh{\widehat{f}}
\def\P{\Cal P}
\def\ek{e^*_k}
\def\eo{e^0_k}
\def\hk{h^*_k}

\def\bbi#1#2{\left[\matrix #1\\#2\endmatrix\right]}
\def\bi#1#2#3{{\bbi {#1} {#2}}_{\!#3}}
\def\bqt#1#2{{\bbi{#1} {#2}}_{\!q,t}}
\def\blm{{\bbi \l \mu}_{\!q,t}}
\def\bl#1{{\bbi l m}_{\!#1}}

\def\bB{\bold B}
\def\bA{\bold A}
\def\bC{\bold C}
\def\bS{\bold S}

\def\li #1 #2 #3 {\vrule height #1 pt width #2 pt depth #3 pt}

\def\hl #1 #2 #3 #4 {\rlap{ \kern #1 pt \raise -#2 pt \hbox{ \li #4 #3 0 }}}

\def\vl #1 #2 #3 #4 {\rlap{ \kern #1 pt \raise -#2 pt \hbox{ \li #3 #4 0 }}}

\def\wr #1 #2 #3 {\rlap{ \kern #1 pt \raise -#2 pt \hbox{\tenrm #3 }}}

\def\picture #1 #2 #3 \endpicture
   { \bigskip \centerline{\hbox to #1 pt
   {\nullfont #3 \hss}} \nobreak \bigskip \centerline{#2} \bigskip}

\def\tht{\thetag}

\topmatter
\title 
Binomial formula for Macdonald polynomials\\
and its applications
\endtitle
\author
Andrei Okounkov
\endauthor
\address
University of Chicago, Dept.\ of Mathematics,
5734 S.\ University Ave., Chicago, IL 60637
\endaddress
\email
okounkov\@math.uchicago.edu
\endemail
\thanks
The paper was written during the author's stay at the
Institute for Advanced Study in Princeton. I am 
grateful to IAS for hospitality and to NSF for
financial support (grant DMS 9304580).
\endthanks
\abstract
We generalize the binomial formula for Jack polynomials 
proved  in \cite{OO2} and consider some applications. 
\endabstract
\endtopmatter

\document

\head
\S1 Binomial formula
\endhead

Binomial type theorems (that is Taylor and Newton
interpolation expansions about various points)
are powerful tools for handling special functions.
The high-school binomial formula is the Taylor expansion of the
function $f(x)=x^l$ about the point $x=1$. Its $q$-deformation
is the Newton interpolation with knots
$$
x=1, q, q^2, \dots\,,
$$
which reads
$$
x^l=\sum_m 
\bl{q}
{(x-1)\dots(x-q^{m-1})} \,, \tag 1.1
$$
where 
$$
\bl{q}=\frac{(q^l-1)\dots(q^l-q^{m-1})}
{(q^m-1)\dots(q^m-q^{m-1})}
$$
is the $q$-binomial coefficient. 
Denote the Newton interpolation polynomials in the RHS of
\tht{1.1} by
$$
P^*_k(x;q)=(x-1)\cdots(x-q^{k-1})\,,\quad k=0,1,2,\dots\,.
$$
The formula inverse to \tht{1.1} is 
the Taylor expansion of $P^*_l(x;q)$ about
the point $x=0$ 
$$
\frac{(x-1)\dots(x-q^{l-1})}{q^{l(l-1)/2}} = 
\sum_m
(-1)^{l-m} 
\bl{1/q}
\frac{x^m}{q^{m(m-1)/2}} \,, \tag 1.2
$$
which is essentially the fundamental 
$q$-binomial theorem \cite{An,GR}.

The formulas \tht{1.1} and \tht{1.2} are limit cases of the formula
$$
\frac{P^*_l(ax;q)} {P^*_l(a;q)}= 
\sum_m
a^m 
\bl{1/q}
\frac{P^*_m(x;1/q)}
{P^*_m(a;q)}  \tag 1.3
$$
as $a\to\infty$ and $a\to 0$ respectively.
Since
$$
\bl{1/q}=
\frac{P^*_m(q^{-l};1/q)}
{P^*_m(q^{-m};1/q)}
$$
the formula \tht{1.3} is symmetric in $q^{-l}$ and $x$.
The formula \tht{1.3} is easy to remember in the form
$$
\frac{P^*_l(ax;q)} {P^*_l(a;q)}= 
\sum_m
\big(\dots\big) \,
P^*_m(q^{-l};1/q) \,
P^*_m(x;1/q)\,, \tag 1.4
$$
where the factor
$$
\big(\dots\big)=
\frac{a^m}
{P^*_m(q^{-m};1/q) \, P^*_m(a;q)} 
$$
can be easily reconstructed using that
\roster
\item this factor does not depend neither on $x$ nor on $l$, and 
\item the highest degree term of $P^*_l(x;q)$ is $x^l$.
\endroster

In this paper we consider a multivariate generalization of 
the formula \tht{1.4} for symmetric functions.
Let $q$ and $t$ be two parameters.
Let $\mu$ be a partition. Recall the following notations of
Macdonald \cite{M}. Set
$$
n(\mu)=\sum_i (i-1) \mu_i = \sum_j \mu'_j(\mu'_j-1)/2 \,.
$$
For each square $s=(i,j)\in\mu$ the numbers
$$
\alignat2
&a(s)=\mu_i-j,&\qquad &a'(s)=j-1,\\
&l(s)=\mu'_j-i,&\qquad &l'(s)=i-1,
\endalignat
$$
are called the arm-length, arm-colength, leg-length, and
leg-colength respectively. We shall need the following
normalization constant
$$
H(\mu;q,t)=t^{-2n(\mu)} q^{n(\mu')} \, \prod_{s\in\mu} 
\left(q^{a(s)+1} t^{l(s)} -1\right)\,. \tag 1.5
$$

Denote by $\P$ the set of all partitions. This set is
partially ordered by inclusion of partitions $\mu\subset\l$.
Consider the following multivariate Newton interpolation
polynomials

\definition{Definition (\cite{Kn,S2,Ok3})} Suppose $\mu$ has $\le n$
parts. The {\it interpolation Macdonald polynomial}
$$
\Pm(x_1,\dots,x_n;q,t)
$$
is the unique polynomial 
in $x_1,\dots,x_n$  such that
\roster
\item $\Pm(x_1,\dots,x_n;q,t)$ is symmetric in variables
$x_it^{-i}$, 
\item $\deg\Pm=|\mu|$,
\item $\Pm(q^\mu;q,t)=H(\mu;q,t)$, and
\item $\Pm(q^\l;q,t)=0$ for all
partitions $\l$ such that $\mu\not\subset\l$. 
\endroster
Here $\Pm(q^\l;q,t)=\Pm(q^{\l_1},\dots,q^{\l_n};q,t)$.
Another name for the same polynomials is the 
{\it shifted Macdonald polynomials}.
\enddefinition

The connection with ordinary Macdonald polynomials will
be explained below. It is easy to see that the polynomial
$\Pm(x;q,t)$ is indeed unique provided it exists (which
is not evident). It follows
from the definition that $\Pm(x;q,t)$ is stable, that is
$$
\Pm(x_1,\dots,x_n,1;q,t)=\Pm(x_1,\dots,x_n;q,t)
$$
provided $\mu$ has at most $n$ parts. Using this stability
one can define polynomials $\Pm(x;q,t)$ in infinitely many
variables.

\remark{Remark} The condition \therosteritem{4} in the 
definition of $\Pm(x;q,t)$ can be replaced by a weaker
condition 
\roster
\item"(4$'$)"
$\Pm(q^\l;q,t)=0\quad\text{for all
partitions $\l\ne\mu$ such that $|\l|\le|\mu|$\,.}$
\endroster
The first proof of the fact that \therosteritem{4$'$}
implies the property \therosteritem{4}, which is in that case
referred to as the {\it extra vanishing}, was given
in \cite{Kn}.
\endremark

The interpolation polynomials in the Schur functions case were studied
in \cite{OO} (and, from a different point of view, in earlier papers 
cited in \cite{OO}). They are called
{\it shifted Schur functions}, have many remarkable properties, and can
be interpreted as distinguished basis elements of the center
of the universal enveloping algebra ${\Cal U}(\frak{gl}(n))$
(see also \cite{Ok1,2} where explicit formulas for the 
corresponding central elements were given). 
We have also conjectured (unpublished) the general
combinatorial formula for the polynomials $\Pm(x;q,t)$ (the 
formula \tht{1.6} below). 
Then F.~Knop and S.~Sahi proved a series of results for
general $(q,t)$ in \cite{KS,Kn,S2}. In particular, it was
proved that the polynomial $\Pm(x;q,t)$ exists 
and has the fundamental property \tht{1.8} below. 
They used difference equations for the polynomials $\Pm(x;q,t)$.
Later the combinatorial formula \tht{1.6} and also
a $q$-integral representation for $\Pm(x;q,t)$
were proved in \cite{Ok3}.

We call a tableau
$T$ on $\mu$ a {\it reverse tableau} if its entries strictly
decrease down the columns and weakly decrease in the rows.
By $T(s)$ denote the entry in the square $s\in\mu$. We have \cite{Ok3}
$$
\Pm(x;q,t)=\sum_T \psi_T(q,t) \prod_{s\in\mu} t^{1-\Ts}
\left(x_{\Ts}-q^{a'(s)} t^{-l'(s)}
\right) \,, \tag 1.6
$$
where the sum is over all reverse tableau on $\mu$ with entries
in $\N$ and $\psi_T(q,t)$ is the same $(q,t)$-weight of a tableau
which enters the combinatorial formula for ordinary Macdonald polynomials
(see \cite{M},\S VI.7)
$$
P_\mu(x;q,t) = \sum_T \psi_T(q,t) \prod_{s\in\mu} x_{\Ts} \,.  \tag 1.7
$$
The coefficients $\psi_T(q,t)$ are rational functions of $q$ and $t$.
Explicit formulas for them are given in \cite{M}, \tht{VI.7.11.${}'$}.
It is clear from \tht{1.6} and \tht{1.7} that
$$
\Pm(x_1,\dots,x_n;q,t)=P_\mu(x_1,x_2 t^{-1},\dots,x_n t^{1-n}) +
\,\text{ lower degree terms} \,. \tag 1.8
$$
In other words, ordinary Macdonald polynomials are the top homogeneous
layer of the interpolation Macdonald polynomials.

It follows from the combinatorial formula \tht{1.6} that
$$
\Pm(a,\dots,a;q,t)= P_\mu(1,t^{-1},\dots,t^{1-n};q,t)
\prod_{s\in\mu} 
\left( a-q^{a'(s)} t^{-l'(s)}
\right)\,, \tag 1.9
$$
where $a$ is repeated $n$ times.  
Recall that  (see \cite{M}, \tht{VI.6.11${}'$})
$$
P_\mu(1,t^{-1},\dots,t^{1-n};q,t)=
t^{n(\mu)+|\mu|(1-n)} \prod_{s\in\mu}
\frac
{1-q^{a'(s)}t^{n-l'(s)}}{1-q^{a(s)}t^{1+l(s)}} \,.
$$
Now we can state our main formula which
is the following generalization of \tht{1.4}

\proclaim{Binomial Theorem}
$$
\frac{\Pl(a x_1,\dots, a x_n;q,t)}
{\Pl(a,\dots, a;q,t)} =
\sum_\mu
\frac{a^{|\mu|}} {t^{(n-1)|\mu|}}
\frac{ \Pm(q^{-\l};1/q,1/t)}
{\Pm(q^{-\mu};1/q,1/t)}
\frac{\Pm(x_n,\dots,x_1;1/q,1/t)}
{\Pm(a,\dots,a;q,t)}\,.
$$ 
\endproclaim

\remark{\bf Remarks}
\roster
\item
By definition of $\Pm(x;q,t)$ only summands with
$\mu\subset\l$ actually enter the sum in the binomial 
theorem.
\item 
By \tht{1.5} and \tht{1.9} all denominators in the binomial theorem have
explicit multiplicative expression.
\item 
One can easily remember the binomial formula in the same 
way as the formula \tht{1.4}.
\item
The binomial formula calls for notation
$$
\blm=\frac{ \Pm(q^{\l};q,t)}
{\Pm(q^{\mu};q,t)} \,.
$$
One can rewrite the binomial formula as an identity for
the above $(q,t)$-binomial coefficients; see 
the formula \tht{2.5} in the next section.
Note that unlike $q$-binomial coefficients the above
$(q,t)$-binomial coefficients  are not polynomials
\footnote
{M.~Lassalle recently pointed out that in particular
case when $\l/\mu$ is a single square these binomial
coefficients appeared in \cite{Ka2}.}.
\item  By \tht{1.9} we have
$$
\qquad
\Pm(a,\dots,a;q,t)=(-a)^{|\mu|} q^{n(\mu')}
t^{-n(\mu)-|\mu|(n-1)} \Pm(1/a,\dots,1/a;1/q,1/t)\,.
$$
Therefore we can rewrite the binomial formula as
follows
$$
\multline
\qquad \frac{\Pl(a x_1,\dots, a x_n;q,t)}
{\Pl(a,\dots, a;q,t)} = \\
\sum_\mu
(-1)^{|\mu|} \frac{t^{n(\mu)}}{q^{n(\mu')}}
\frac{ \Pm(q^{-\l};1/q,1/t)}
{\Pm(q^{-\mu};1/q,1/t)}
\frac{\Pm(x_n,\dots,x_1;1/q,1/t)}
{\Pm(1/a,\dots,1/a;1/q,1/t)}\,. 
\endmultline \tag 1.10
$$ 
\endroster
\endremark

Using \tht{1.8} and the well known equality
$$
P_\l(x;q,t)=P_\l(x;1/q,1/t)
$$
we obtain as limit cases as $a\to\infty$ or
$a\to 0$  binomial formulas involving ordinary
Macdonald polynomials. We have
$$
\frac{P_\l(x_1,\dots,x_n t^{1-n};q,t)}
{P_\l(1,\dots,t^{1-n};q,t)} =
\sum_\mu
\frac{ \Pm(q^{\l};q,t)}
{\Pm(q^{\mu};q,t)}
\frac{\Pm(x_1,\dots,x_n;q,t)}
{P_\mu(1,\dots,t^{1-n};q,t)}\,. \tag 1.11
$$ 
and
$$
\frac{\Pl(x_1,\dots, x_n;q,t)}
{\Pl(0,\dots,0;q,t)} =
\sum_\mu
\frac{ \Pm(q^{-\l};1/q,1/t)}
{\Pm(q^{-\mu};1/q,1/t)}
\frac{P_\mu(x_1,\dots,x_n t^{1-n};q,t)}
{\Pm(0,\dots,0;q,t)}\,, \tag 1.12
$$ 
as $a\to\infty$ and $a\to0$ respectively.

If we set $t=q^\theta$ and let $q\to 1$ then \tht{1.11}
turns into the binomial formula for Jack polynomials
proved in \cite{OO2}. Thus, the above binomial theorem is a
(far reaching) generalization of that formula. 
Various results on binomial coefficients for Jack
polynomials were obtained in \cite{Bi,FK,Ka,Lasc,La,OO}.

In the next section we list some 
applications of the binomial formula. 
Those applications are quite different from the 
applications of the Jack degeneration of the binomial
formula considered in \cite{OO2,OO4,Ok5,Ok6} (see also
\cite{KOO}). Then in section 3
we recall some results of F.~Knop and S.~Sahi. The proof
of the theorem is given in section 4. Finally, in Appendix,
we discuss an elementary algorithm for multivariate
symmetric Newton interpolation (to avoid possible
confusion: this is not the interpolation considered
by Lascoux and Sch\"utzenberger.). 

A binomial formula for the Koornwinder polynomials
(which generalize Macdonald polynomials for the
classical root systems) was proved in \cite{Ok4}. 

I would like to thank G.~Olshanski and S.~Sahi for
helpful discussions.

\head
\S2 Applications
\endhead

Consider the following matrix whose rows and columns
are indexed by partitions with at most $n$ parts
$$
\bS(a,n;q,t)=
\left(\,
\frac{\Pl(a q^{-\nu_n},\dots, a q^{-\nu_1};q,t)}
{\Pl(a,\dots, a;q,t)}
\,\right)_{\l,\nu} \,.
$$
{F}rom the binomial theorem we have
\proclaim{Symmetry} The matrix $\bS(a,n;q,t)$ is symmetric
for any $a$.
\endproclaim

This symmetry is a simple yet central result. Since 
it interchanges the label and the argument of the polynomial
$\Pm$ all properties of $\Pm$ come in pairs, where
one is obtained from the other using the above symmetry.
For example, one checks that we have two following pairs
$$
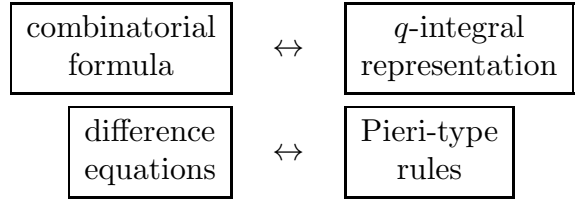
\align
\boxed{
\matrix
\text{combinatorial}\\
\text{formula}
\endmatrix} \quad
&\leftrightarrow
\quad
\boxed{
\matrix
\text{$q$-integral}\\
\text{representation}
\endmatrix} \\
\boxed{
\matrix
\text{difference}\\
\text{equations}
\endmatrix} \quad
&\leftrightarrow
\quad
\boxed{
\matrix
\text{Pieri-type}\\
\text{rules}
\endmatrix}
\endalign
$$

Similarly, it is evident from \tht{1.11} that the matrix
$$
\bS(\infty,n;1/q,1/t)=
\left(\,
\frac{P_\l( q^{\nu_1},\dots,q^{\nu_n}t^{1-n};q,t)}
{P_\l(1,\dots,t^{1-n};q,t)}
\,\right)_{\l,\nu} \tag 2.1
$$
is symmetric, which was conjectured by I.~Macdonald, 
proved by T.~Koornwinder (unpublished, see \cite{M}, section
VI.6), and, by different methods,
in \cite{EK} and \cite{Ch}.

The binomial formula is just the Gauss decomposition
of the symmetric matrix 
$$
\bS(a,n;q,t) = \bB'(1/q,1/t) \, \bA(a,n;1/q,1/t) \, \bB (1/q,1/t)
\tag 2.2
$$
where prime stands for transposition, the matrix
$$
\bB(q,t)=\left( \blm
\right)_{\mu,\l} \tag 2.3
$$
is upper triangular (with respect to the partial ordering of
partitions by inclusion) and unipotent, and the matrix
$$
\bA(a,n;q,t) = \left( \delta_{\l\mu} 
(-1)^{|\mu|} \frac{q^{n(\mu')}}{t^{n(\mu)}}
\frac
{\Pm(q^{\mu};q,t)}
{\Pm(1/a,\dots,1/a;q,t)}
\,\right)_{\mu,\l} 
$$
is diagonal.

The symmetry of the matrix $\bS$ has the following
combinatorial interpretation. Suppose that for some $b\in\N$
$$
\l,\nu \subset \beta
$$
where $\beta$ is a box $\beta=(b^n)$. 
Denote by
$$
\beta\setminus\lambda=(b-\l_n,\dots,b-\l_1)
$$
the partition complementary to $\l$ in $\beta$. Setting $a=q^b$
we obtain 
$$
\frac{\Pl(q^{\beta\setminus\nu};q,t)}
{\Pl(q^{\beta};q,t)} =
\frac{P^*_\nu(q^{\beta\setminus\l};q,t)}
{P^*_\nu(q^{\beta};q,t)} \,. \tag 2.4
$$
The skew diagram $(\beta\setminus\nu)/\lambda$ is
depicted in the following Fig.~1 (where the
diagram $\nu$ is rotated by $180^\circ$)
\picture 125 {Fig.~1}
\hl 0 0 125 0.15
\hl 70 5 55 0.15
\hl 60 10 65 0.15
\hl 55 15 70 0.15
\hl 45 20 75 0.15
\hl 35 25 85 0.15
\hl 20 30 95 0.15
\hl 15 35 100 0.15
\hl 15 40 95 0.15
\hl 10 45 90 0.15
\hl 5 50 90 0.15
\hl 5 55 90 0.15
\hl 0 60 85 0.15
\hl 0 65 75 0.15
\hl 0 70 55 0.15
\hl 0 75 125 0.15
\vl 0 75 75 0.15 
\vl 5 75 25 0.15
\vl 10 75 30 0.15
\vl 15 75 40 0.15
\vl 20 75 45 0.15
\vl 25 75 45 0.15
\vl 30 75 45 0.15
\vl 35 75 50 0.15
\vl 40 70 45 0.15
\vl 45 70 50 0.15
\vl 50 70 50 0.15
\vl 55 70 55 0.15
\vl 60 65 55 0.15
\vl 65 65 55 0.15
\vl 70 65 60 0.15
\vl 75 65 60 0.15
\vl 80 60 60 0.15
\vl 85 60 60 0.15
\vl 90 55 55 0.15
\vl 95 55 55 0.15
\vl 100 45 45 0.15
\vl 105 40 40 0.15
\vl 110 40 40 0.15
\vl 115 35 35 0.15
\vl 120 25 25 0.15
\vl 125 75 75 0.15
\wr 15 17 {$\lambda$}
\wr 110 64 {$\nu$}
\wr 130 75 {$\beta$}
\endpicture
\noindent 
Introduce the following {\it trinomial coefficients}
$$
\bqt{\beta}{\l,\nu} =
\bqt{\beta}{\nu}
\bqt{\beta\setminus\nu}{\l}\,.
$$
Using them one rewrites \tht{2.4} as follows
\proclaim{Symmetry II}
$$
\bqt{\beta}{\l,\nu} =
\bqt{\beta}{\nu,\l} \,.
$$
\endproclaim

In particular,
$$
\bqt{\beta}{\nu}=
\bqt{\beta}{\beta\setminus\nu}\,.
$$
In the shifted Schur function case the symmetry
is obvious (see Fig.~1) 
from the following formula (see \cite{OO}, \S8, and also
\cite{OO2}, \S5 for Jack generalization)
$$
\bi{\beta}{\l,\nu}{1,1}= 
\binom {|\beta|}{|\l|,|\nu|}\, \frac{\dim\l\, \dim\nu}{\dim \beta}
\, \dim (\beta\setminus\nu)/\lambda\,,
$$
where $\dim$ stands for the number of standard tableaux
on a (skew) diagram.

Observe that 
$$
\bA(q^b,n;1/q,1/t)=
\left( \delta_{\l\mu} 
(-1)^{|\mu|} \frac{t^{n(\mu)}}{q^{n(\mu')}}
\bi{\beta}{\mu}{1/q,1/t}^{-1}
\,\right)_{\mu,\l}
$$
Therefore, the binomial formula can be rewritten as follows

\proclaim{Binomial Theorem II}
$$
\bqt{\beta}{\l,\nu}=
\bqt{\beta}{\l} \bqt{\beta}{\nu} 
\sum_\mu (-1)^{|\mu|} 
\frac{t^{n(\mu)}}{q^{n(\mu')}}
\bi{\l}{\mu}{1/q,1/t}
\bi{\nu}{\mu}{1/q,1/t}
\bi{\beta}{\mu}{1/q,1/t}^{-1} \,. \tag 2.5
$$
\endproclaim
Note that only summands with $\mu\subset\l\cap\nu$
actually enter this sum.

We denote by $\Ltn$ the algebra of polynomials in $x_1,\dots,x_n$
which are symmetric in variables 
$$
x_i t^{-i} \,.
$$
We call such polynomials {\it shifted symmetric}. 
The polynomials $\Pm(x;q,t)$, where $\mu$ ranges over partitions
with at most $n$ parts, form a linear basis in this algebra.

The binomial formula gives transition coefficients between
various mutually triangular linear bases of $\Ltn$.
Combining, for example, \tht{1.11}
and \tht{1.12} and taking into account that by \tht{1.9}
$$
\frac{\Pm(0,\dots,0;q,t)}
{P_\mu(1,t^{-1},\dots,t^{1-n};q,t)}=
(-1)^{|\mu|} q^{n(\mu')} t^{-n(\mu)}
$$
we obtain
$$
\delta_{\l\nu}=(-1)^{|\nu|}  \frac{t^{n(\nu)}}{q^{n(\nu')}}
\sum_{\mu} (-1)^{|\mu|} \frac{q^{n(\mu')}}{t^{n(\mu)}}
\bi{\l}{\mu}{q,t}\bi{\mu}{\nu}{1/q,1/t}\,.
$$
This is equivalent to the following proposition

\proclaim{Inversion} The inverse matrix of the matrix \tht{2.3}
is given by
$$
\bB^{-1}(q,t)= \bC^{-1}(q,t)\, \bB(1/q,1/t)\, \bC(q,t)\,, \tag 2.6
$$
where $\bC(q,t)$ is the diagonal matrix
$$
\bC(q,t)=
\left( \delta_{\l\mu} (-1)^{|\mu|} \,
\frac{q^{n(\mu')}}{t^{n(\mu)}} \,
\right)_{\mu,\l}\,.
$$
\endproclaim

Note that the inverse of the (truncated) matrix \tht{2.1} in the case
when $q$ is a root of unity was calculated by A.~A.~Kirillov, Jr.\ 
\cite{Ki} and also by I.~Cherednik \cite{Ch}.

The algebra $\Ltn$ is filtered by degree of polynomials. We
denote by $\Lt$ the projective limit of algebras $\Ltn$, $n\in\N$ 
with respect to homomorphisms
$$
f(x_1,\dots,x_{n+1}) \mapsto f(x_1,\dots,x_n,1)\,.
$$
Note that $f(q^\l)$ is well defined for all $f\in\Lt$ and $\l\in\P$.
Consider the Newton interpolation of a polynomial $f\in\Lt$ with
knots 
$$
x=q^\l\,,\quad \l\in\P\,,
$$
or, in other words, the expansion of $f$ in
the linear basis $\{\Pm(x;q,t)\}_{\mu\in\P}$ of $\Lt$.
The coefficients $\fh(\mu;q)$ of this expansion
$$
f(x)=\sum_\mu \fh(\mu;q)\, \Pm(x;q,t)
$$
can be found  from the following non-degenerate
triangular (with respect to ordering of $\P$ by inclusion)
system of linear equations
$$
f(q^\l)=\sum_{\mu\subset\l} \fh(\mu;q)\, \Pm(q^\l;q,t)\,,
\quad \mu,\l\in\P\,. \tag 2.7
$$

Using \tht{2.6} we can solve the equations \tht{2.7} and obtain 

\proclaim{Explicit Newton Interpolation}
The coefficients $\fh(\mu;q)$ are given by
$$
\fh(\mu;q)=\frac1{\Pm(q^\mu;q,t)}
\sum_{\nu\subset\mu}
(-1)^{|\mu/\nu|}
\frac{q^{n(\mu'/\nu')}}
{t^{n(\mu/\nu)}}
\bi{\mu}{\nu}{1/q,1/t}\, f(q^\nu) \,, \tag 2.8
$$
where $n(\mu/\nu)=n(\mu)-n(\nu)$.
\endproclaim

In particular, this formula can be used for expansion 
of any homogeneous symmetric polynomial in ordinary Macdonald
polynomials.

One can think of the function 
$$
\fh: \mu \mapsto \fh(\mu,q)\, \quad \mu\in\P
$$
is an of analog of the {\it Fourier transform} of the function
$$
f:  \mu \mapsto f(q^\mu)\, \quad \mu\in\P \,.
$$
The inverse transform is obviously given by \tht{2.7}.

Combining the formula \tht{2.8} with the Knop-Sahi difference
equation for the polynomials $\Pm$ one can obtain a nice
multivariate generalization of the classical {\it algorithm}
for Newton interpolation with knots $0,1,2,\dots$ 
(see Appendix).

The binomial theorem provides a natural proof
of the {\bf $q$-integral representation} for interpolation Macdonald
polynomials found in \cite{Ok3}. In \cite{Ok3} the integral
representation was guessed and checked.
One can {deduce} the integral representation from
the formulas
\tht{1.7} and \tht{1.11} in exactly the same way
as it was done in \cite{OO2} for Jack polynomials.  

Notice that the the binomial formula makes certain 
sense for $\l$ which is  not a partition. An example of such 
{\it analytic continuation} is the formula \tht{2.10} below. 
Notice, however, that the RHS may or may not converge to the 
LHS in the usual analytical sense.

As two particular cases of the binomial formula, one obtains two
following generation functions for the 1-column and 1-row interpolation
Macdonald polynomials. By definition, set:
$$
\align
\ek(x;t)&=P^*_{(1^k)}(x;q,t)\,, \\
\hk(x;q,t)&=P^*_{(k)}(x;q,t)\,,
\quad k=1,2,\dots \,.
\endalign
$$
The polynomials $\ek(x;t)$ do not depend on $q$ for it follows from 
the combinatorial formula \tht{1.6} that
$$
\ek(x;t)=\sum_{i_1<\dots<i_k} t^{k-\sum i_s} \prod_{s} (x_i-t^{s-k})\,.
$$
Let $u$ be a parameter. We have the following

\proclaim{Generating functions} 
$$
\align 
&\prod_{i=1}^{\infty} 
\frac{1+x_i t^{1-i}/u}
{1+ t^{1-i}/u}
=\sum_{k=0}^{\infty} 
\frac{\ek(x;t)}
{(u+1)(u+1/t)\dots(u+1/t^{k-1})}\,, \tag 2.9 \\
&\prod_{i=1}^{\infty} 
\frac{(x_i t^{1-i}/u;q)_\infty}
{(x_i t^{-i\hphantom{1}}/u;q)_\infty}
\frac
{(t^{-i\hphantom{1}}/u;q)_\infty}
{(t^{1-i}/u;q)_\infty}
=\sum_{k=0}^{\infty} t^{-k}
\frac{(t;q)_k}{(q;q)_k}
\frac{\hk(x;q,t)}
{(u-1)(u-q)\dots(u-q^{k-1})}\,.  \tag 2.10
\endalign
$$
\endproclaim

Here $(z;q)_r=(1-z)(1-qz)\dots(1-q^{r-1}z)$. 
Note that above expansions make sense as formal power series
in $1/u$ with coefficients in $\Lt$.

This generating function \tht{2.9}
is a minor modification of the generating function 
in \cite{OO}, theorem 12.1. 
The analog of the  generating function \tht{2.10} for shifted Jack
polynomials was found by G.~Olshanski (unpublished) and it
was used by him to prove the Jack degeneration of the
1-row case of the combinatorial formula \tht{1.6}.

It is clear that if
$$
x_{n+1}=x_{n+2}=\dots=1
$$
then \tht{2.9} is just the binomial formula for 
$\l=(1^n)$ and $a=-t^{1-n}/u$, whereas \tht{2.10} is the binomial
expansion of
$$
\frac{\Pl(x_n/u,\dots,x_1/u;1/q,1/t)}
{\Pl(1/u,\dots,1/u;1/q,1/t)}
$$
analytically continued to the point $\l$ satisfying
$$
q^\l=(1/t,\dots,1/t)\,.
$$
It is easy to check both expansions directly.
\footnote{ In fact, one can generalize the expansion \tht{2.9}
by replacing the sequence $1,1/t,1/t^2,\dots$ by an arbitrary 
sequence. That type of expansions was considered by 
A.~Abderrezzak in
\cite{Ab}. However, only very few facts about interpolation Macdonald
polynomials admit such a vast generalization, see \cite{Ok7}.}
Since it will be convenient for us to use \tht{2.9} below
in section  3, we shall give such a direct proof, which is
essentially the argument from \cite{OO}, section \S 12. Same 
argument applies to \tht{2.10} 

\demo{Proof of \tht{2.9}} Let $\eo(x;t)$ be the coefficients in the expansion
$$
\prod_{i=1}^{\infty} 
\frac{1+x_i t^{1-i}/u}
{1+ t^{1-i}/u}
=\sum_{k=0}^{\infty} 
\frac{\eo(x;t)}
{(u+1)(u+1/t)\dots(u+1/t^{k-1})}\,. \tag 2.11
$$
We want to prove that $\eo(x;t)=\ek(x;t)$.
It is clear that $\eo\in\Lt$. Set $x_i=\xi_i u z$,
where $\xi_i$ and $z$ are new variables and let $u\to\infty$.
One obtains
$$
\eo(x;t)=e_k(x_1,x_2 t^{-1},\dots) +
\,\text{ lower degree terms}\,. 
$$
Therefore, it suffices to prove that
$$
e^0_n(x_1,\dots,x_{n-1},1,\dots,1;t)=0\,,\quad n=1,2,\dots\,.
$$
Substitute $x_{n}=x_{n+1}=\dots=1$ in \tht{2.11}. Then the
RHS becomes
$$
\prod_{i=1}^{n-1} 
\frac{u+x_i t^{1-i}}
{u+t^{1-i}}\,.
$$
In its expansion of the form \tht{2.11} all summands with $k\ge n$
vanish. \qed
\enddemo

The expansions \tht{2.9} and \tht{2.10} are related by the {\bf duality}
for the interpolation Macdonald polynomials,
see \cite{Ok3}, theorem IV. That duality (in particular, the proposition 
6.1 from \cite{Ok3})
implies the following duality for $(q,t)$-binomial
coefficients
$$
\bi{\l}{\mu}{q,t}=
\bi{\l'}{\mu'}{1/t,1/q}\,. \tag 2.12
$$

\head
\S3 Difference equations for $\Pm(x;q,t)$.
\endhead

In this section we recall some results due to F.~Knop and S.~Sahi
\cite{Kn,S2}.
Only the explicit formulas for higher order difference operators are new.

Denote by $T_{x,q}$ the shift operator
$$
[T_{x,q} f](x)=f(q x)\,.
$$
Observe that for any polynomials $p_1,p_2$ in $x$ such that
$$
p_1(0)=p_2(0)
$$
the operator
$$
f\longmapsto \frac{ p_1(x) f(x)- p_2(x) f(q x)}{x}
$$
maps polynomials to polynomials. Set
$$
T_i = T_{x_i,1/q}\,.
$$
Consider the following $n$-dimensional difference operator 
$$
D(u)=  \Vm \det \left[ 
\frac{(1-u x_i)}{x_i} (1-t x_i)^{j-1} - 
\frac{(1-x_i)^j}{x_i} \, T_i \,
\right]_{1\le i,j\le n} \,. 
$$
Here $u$ is a parameter and $V(x)$ is the Vandermonde determinant.
The determinant is well defined since multiplication
by $x_i$ and the shift $T_j$ commute provided $i\ne j$. It is
clear that $D(u)$ maps $\Ln$ to $\Ln$ and does not increase degree.

For any subset $I\subset\{1,\dots,n\}$ set
$$
T_I = \prod_{i\in I} T_i \,.
$$
Using the Vandermonde determinant formula one expands $D(u)$
as follows
$$
D(u)=\xm \sum_{I\subset\{1,\dots,n\}} (-1)^{|I|} t^{\dots}
\prod_{i\in I} (1-x_i)
\prod_{i\notin I} (1-u x_i)
\prod_{i\in I,\,j \notin I}
\frac{x_i-t x_j} {x_i - x_j}\,\, T_I \,\,, 
$$
where the omitted exponent of $t$ equals $(n-|I|)(n-|I|-1)/2$.

Now shift the variables
$$
x_i\mapsto x_i t^{n-i}
$$
and consider the operator
$$
\multline
D^*(u)=\xm \sum_{I\subset\{1,\dots,n\}} (-1)^{|I|}\, 
t^{(n-|I|)(n-|I|-1)/2} \\
\prod_{i\in I} (1-x_i t^{n-i})
\prod_{i\notin I} (1-u x_i t^{n-i})
\prod_{i\in I,\,j \notin I}
\frac{x_i - x_j t^{i-j+1}} 
{x_i - x_j t^{i-j}}\,\, T_I \,\,, 
\endmultline 
$$
which maps $\Ltn$ to $\Ltn$.

Below we shall have to evaluate
$$
\left[
D^*(u) f
\right](x)\,, \quad f\in \Ltn
$$
at the points of the form
$$
x=(a q^{\xi_1},\dots,a q^{\xi_n})\,,\quad \xi_i\in \Z\,, \tag 3.1
$$
where $\xi_1\ge\dots\ge\xi_n$ and $a$ is a parameter (which will
be the same parameter $a$ 
as in the binomial theorem). Fix a subset $I$ and set
$$
\xi'_i=
\cases
\xi_i -1, & i\in I\\
\xi_i, & i\notin I\,.
\endcases
$$
We have the following obvious 

\proclaim{Lemma} Let $x$ be as in \tht{3.1}. Then
$$
\prod_{i\in I,\,j \notin I}
\frac{x_i - x_j t^{i-j+1}} 
{x_i - x_j t^{i-j}} \ne 0
$$
iff $\xi'_1\ge\dots\ge\xi'_n$
\endproclaim

This lemma results in the following

\proclaim{Theorem}
$$
D^*(u)\, \Pm(x;q,t) = 
\left( \prod_i (q^{-\mu_i} t^{i-1} - u t^{n-1})
\right) \, \Pm(x;q,t)\,. \tag 3.2
$$
\endproclaim

This theorem (without explicit formula for $D^*(u)$) is due
to F.~Knop and S.~Sahi \cite{Kn,S2}.

\demo{Idea of Proof} 
Using that
$$
\prod_{i\in I} (1-x_i t^{n-i}) = 0 
\quad \text{if} \quad x_n=1 
\quad \text{and} \quad i\in I
$$
and the lemma one easily checks that
$$
\left[ D^*(u) \Pm
\right] (q^\l)= 0,\quad \l\in\P,\, \l\not\subset\mu
$$
and that
$$
\left[ D^*(u) \Pm
\right] (q^\mu) = 
\left( \prod_i (q^{-\mu_i} t^{i-1} - u t^{n-1}) 
\right)\Pm(q^\mu)  \,. \qed
$$
\enddemo

In particular, the operators $D^*(u)$ for different
values of $u$ commute. Consider the algebra $\D$
generated by these commuting difference operators.

This commutative algebra is
isomorphic to the algebra $\Lotn$ under the Harish-Chandra
map which maps an difference operator $D\in\D$ to
the polynomial $d \in \Lotn$ such that 
$$
D \,\Pm(x;q,t) = d(q^{-\mu}) \, \Pm(x;q,t) \,.
$$

Let $D_k$ be the coefficient 
$$
D^*(u) = \sum_{k\le n} (-1)^{n-k} t^{(n-1)(n-k)}
(u-1)(u-1/t)\dots(u-1/t^{n-k-1}) \, D_k
$$
in the Newton interpolation 
of $D^*(u)$ with knots 
$$
u=1,1/t,1/t^2,\dots\,.
$$
By \tht{3.2} and \tht{2.9} we have
$$
D_k \,\Pm(x;q,t) = e^*_k(q^{-\mu};1/t) \, \Pm(x;q,t) \,.
$$
The operators $D_k$ are generators of the algebra $\D$.

\head
\S4 Proof of the binomial theorem.
\endhead

We begin with general remark about the proof.
As it is usual in the Macdonald polynomials theory,
various pieces of the puzzle can be arranged in
various combination to give a proof. Below
we shall use the approach based on the Knop-Sahi
difference equations. 
However, it will be  clear that the formula \tht{4.4},
which states that the LHS is a polynomial in $q^{-\l}$
with a certain type of symmetry, is the key point
of the proof. Instead of using difference
equations, one  can deduce \tht{4.4} and hence
the binomial theorem from the $q$-integral representation
of interpolation Macdonald polynomials \cite{Ok3}.  

In this section by $c_\ga$ we always denote {\it some} coefficients.
This notation is used to show which summands possibly enter the formula.

It is clear from definitions that
$$
D_k=
\sum_{|I|=k} (-1)^{|I|} 
t^{\dots} 
\prod_{i\in I} \frac{(1-x_i t^{n-i})}{x_i}
\prod_{i\in I,\,j \notin I}
\frac{x_i - x_j t^{i-j+1}} 
{x_i - x_j t^{i-j}}\,\, T_I \, +
\sum_{|I|<k}\dots\,. 
$$
Suppose
$$
x=(a q^{\xi_1},\dots,a q^{\xi_n})\,,\quad \xi\in \Z\,,
$$
where $\xi_1\ge\dots\ge\xi_n$. Then by the lemma
$$
\left[
D_k f \right] (x) =
\sum_{\eta} c_\eta f(a q^{\eta_1},\dots,a q^{\eta_n}) \tag 4.1
$$
with some coefficients $c_\eta$, where 
$$
\align
& \eta_1\ge \dots \ge \eta_n\,,\\
& \xi_i - \eta_i \in \{0,1\}\,, \quad i=1,\dots,n\,,\tag 4.2\\
& \sum \xi - \sum \eta_i \le k \,.
\endalign
$$
Order the $n$-tuples $\eta$ lexicographically as follows
$$
\eta > \te
$$
if $\eta_n < \te_n$, or if $\eta_n = \te_n$ and
$\eta_{n-1} < \te_{n-1}$, and so on. Then
the maximal $\eta$ satisfying \tht{4.2} is
$$
\etm=(\xi_1,\dots,\xi_{n-k},\xi_{n-k+1}-1,\dots,\xi_{n}-1)\,.
$$
It is clear that the coefficient 
$$
c_{\etm}\ne 0\
$$
\
in \tht{4.1} is a non-zero rational function of $a$.

Now let $\mu$ be a partition and let $\mu'$ be the conjugate
partition. Iterating the above argument we see that
$$
\left[
\prod_i D_{\mu_i} f 
\right] (a,\dots,a) =
\sum_{\nu} c_\nu f(a q^{-\nu_n},\dots, a q^{-\nu_1})\,, \tag 4.3
$$
where $c_\nu$ are some coefficients and $\nu$ ranges over
partitions such that $|\nu|\le|\mu|$ and 
$$
\nu \le \mu'
$$
in the lexicographic order of partitions. Observe that the
coefficient $c_{\mu'}$ in \tht{4.3} does not vanish. Therefore
there exist an operator $D(\mu,a)\in\D$
$$
D(\mu,a)=\sum_{|\nu|\le|\mu|,\nu\le\mu} c_\nu 
\prod_i D_{\nu'_i} 
$$
such that
$$
\left[
D(\mu,a) f 
\right] (a,\dots,a) =
f(a q^{-\mu_n},\dots, a q^{-\mu_1})\,. 
$$
Denote by $d_{\mu,a}\in\Lotn$ the Harish-Chandra image of $D(\mu,a)$
$$
D(\mu,a)\, \Pl(x;q,t) = d_{\mu,a}(q^{-\l}) \, \Pl(x;q,t) \,.
$$
By definition of $D(\mu,a)$ we have
$$
\frac{\Pl(a q^{-\mu_n},\dots, a q^{-\mu_1})}
{\Pl(a,\dots, a)} = d_{\mu,a}(q^{-\l})
\,. \tag 4.4
$$
Since the degree of $e^*_k$ is $k$ we have
$$
\deg d_{\mu,a} \le |\mu|\,.
$$

Now consider the Newton interpolation of the function
$$
f(x)=\frac{\Pl(a x_1,\dots, a x_n;q,t)}
{\Pl(a,\dots,a;q,t)} \quad \in \Ltn
$$
with knots
$$
x=(q^{-\mu_n},\dots, q^{-\mu_1}), \quad \mu\in\P
$$
This interpolation
has the form
$$
\frac{\Pl(a x_1,\dots, a x_n;q,t)}
{\Pl(a,\dots, a;q,t)} =
\sum_{\mu} b_{\mu,\l,a} \, \Pm(x_n,\dots,x_1;1/q,1/t) \,. \tag 4.5
$$
The coefficients $b_{\mu,\l,a}$ are linear combinations of
the values
$$
\frac{\Pl(a q^{-\nu_n},\dots, a q^{-\nu_1};q,t)}
{\Pl(a,\dots, a;q,t)}\,, \quad \nu\subset\mu\,.
$$
By \tht{4.4} we have
$$
b_{\mu,\l,a} = b_{\mu,a} (q^{-\l})\,,
$$
for certain polynomials
$$
b_{\mu,a}\in\Lotn\,,\quad \deg b_{\mu,a} \le |\mu|\,.
$$

Next observe that the highest degree term
of the LHS of \tht{4.5} equals
$$
\frac{a^{|\l|}}
{\Pl(a,\dots,a;q,t)}
P_\l(x_1,x_2 t^{-1},\dots, x_n t^{1-n};q,t)
$$
and 
$$
\Pm(x_n,\dots,x_1;1/q,1/t) = t^{(n-1)|\mu|}
P_\mu(x_1,x_2 t^{-1},\dots, x_n t^{1-n};q,t) + \dots\,,
$$
where dots stand for lower degree terms. Therefore
$$
b_{\mu,a} (q^{-\l})=
\cases
0\,,&|\lambda|\le|\mu| \quad\text{and}\quad \mu\ne\lambda\,,\\
{a^{|\mu|} t^{(1-n)|\mu|}}
{P^*_\mu(a,\dots,a;q,t)}^{-1}\,, &\mu=\lambda\,.
\endcases
$$
Since $\deg b_{\mu,a}\le|\mu|$ the polynomial $b_{\mu,a}$
is uniquely determined by its values at the points $x=q^{-\l}$,
$|\l|\le|\mu|$. This implies
$$
b_{\mu,a} (x) =
\frac{a^{|\mu|} t^{(1-n)|\mu|} }
{\Pm(a,\dots,a;q,t) }
\frac{ \Pm(x;1/q,1/t)}
{\Pm(q^{-\mu};1/q,1/t)}
$$
and concludes the proof of the theorem.

\def\bs{\bigstar}
\def\o{\omega}
\def\nea{{\tsize \nearrow}}

\def\th{\theta}

\head
Appendix. Symmetric Newton's algorithm
\endhead

The formula \tht{2.8} can be turned into a quite
efficient algorithm for multivariate symmetric Newton
interpolation.

Reccurence relations for generalized binomial coefficients
were first used in the Jack polynomial case by M.~Lassalle,
see \cite{La}
\footnote{See also M.~Lassalle's recent preprint
{\it Coefficients binomiaux g\'en\'eralis\'es et
polyn\^omes de Macdonald}, \S10.}. Those relations correspond 
to a Pieri type formula for interpolation Macdonald
polynomials (see \cite{OO2}, \S5). They are related
to the Knop-Sahi type relations used below by
the {\it symmetry}, see \S2.

Let $\l$ be a diagram and suppose that $\nu$
is obtained from $\l$ by removing a corner. 
We shall denote this by writing
$$
\nu\nea\l\,.
$$
Denote the corner being removed  by $\bs$, see Fig.~2
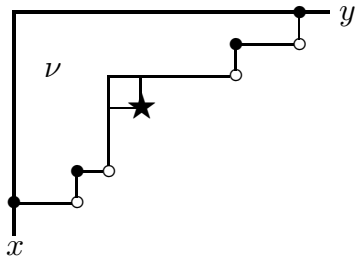
\picture 100 {Fig.~2.\ The diagram $\nu$ is obtained from $\l$ by removing the corner $\bs$}
\hl 0 0 100 1
\vl 0 70 70 1
\hl 0 60 18.5 0.5
\vl 20 58 8 0.5
\hl 20 50 8.5 0.5 
\vl 30 48 28 0.5
\hl 30 20 38.5 0.5 
\vl 70 18 8 0.5
\hl 70 10 18.5 0.5 
\vl 90 8 8 0.5
\hl 30 30 10 0.5 
\vl 40 30 10 0.5
\wr 36 32 {$\bigstar$}
\wr 18 52 {$\bullet$}
\wr -2 62 {$\bullet$}
\wr 68 12 {$\bullet$}
\wr 88 2 {$\bullet$}
\wr 28 52 {$\circ$}
\wr 18 62  {$\circ$}
\wr 68 22 {$\circ$}
\wr 88 12 {$\circ$}
\wr -2 76 {$x$}
\wr 103 2 {$y$}
\wr 10 20 {$\nu$}
\endpicture
\noindent
Consider the outer corners of the diagram $\nu$
(denoted by $\circ$ in Fig.~2) and also the inner
corner of $\nu$ excluding the one where the 
square $\l/\nu$ is attached (those inner
corners are denoted by $\bullet$ in Fig.~2).
By definition, set 
$$
\o(\nu\nea\l;q,t)=
\frac{
\dsize \prod_{\text{inner corners $\bullet$ \,}}
\left(1-q^{\Delta y(\bigstar,\bullet)} t^{\Delta x(\bullet, \bigstar)}\right)}
{
\dsize \prod_{\text{outer corners $\circ$}}
\left(1-q^{\Delta y(\bigstar,\circ)} t^{\Delta x(\circ, \bigstar)}\right)}\,,
$$
where 
$$
\align
\Delta y(\bigstar,\bullet)&=\text{$y$-coordinate $(\,\bigstar\,)$} -
\text{$y$-coordinate $(\,\bullet\,)$}\,,\\
\Delta x(\bullet,\bigstar)&=\text{$x$-coordinate $(\,\bullet\,)$} -
\text{$x$-coordinate $(\,\bigstar\,)$}\,.
\endalign
$$

Now let $f\in\Lt$ be the shifted symmetric polynomial to be 
interpolated. For any pair $\mu\subset\l$ set by definition
$$
F(\mu,\l)=
(-1)^{|\l/\mu|}
\frac{q^{n(\l'/\mu')}}
{t^{n(\l/\mu)}}
\bi{\l}{\mu}{1/q,1/t}\, f(q^\mu) \,,
$$
where $n(\l/\mu)=n(\l)-n(\mu)$.
Then, the interpolation formula \tht{2.8} can be 
rewritten as follows
$$
\fh(\l)=\sum_{\mu\subset\l} F(\mu,\l) \,. \tag 5.1
$$

On easily computes that the Knop-Sahi difference equation
$$
D_1 P^*_\mu(x;q,t)= \left(\sum_i t^{i-1} (q^{-\mu_i}-1)
\right) P^*_\mu(x;q,t)
$$
implies that
$$
\sum_{\nu\nea\l} \o(\nu\nea\l;q,t) \bi{\nu}{\mu}{1/q,1/t} =
\left(\bi{\l}{1}{q,t} - \bi{\mu}{1}{q,t}\right) 
\bi{\l}{\mu}{1/q,1/t} \,. \tag 5.2
$$
Therefore, the function $F(\mu,\l)$ satisfies 
the following reccurence relation
$$
F(\mu,\l)=-\sum_{\mu\subset\nu\nea\l} 
\frac{q^{a'(\l/\nu)}}{t^{l'(\l/\nu)}}
\,
\frac{\o(\nu\nea\l;q,t)}
{\tsize \bi{\l}{1}{q,t} - \bi{\mu}{1}{q,t}
} F(\mu,\nu)\,,
$$
which, together with the obvious initial condition
$$
F(\mu,\mu)=f(q^\mu),
$$
give a simple algorithm for calculation of this function.
The summation \tht{5.1} can be performed simultaneously
with the computation of $F(\mu,\l)$. 

Both in formula \tht{2.8} and in the above algorithm
one has to run a cycle over all pairs of partitions
$\mu\subset\l$. This means of order $N^2$ iterations
to interpolate in $N$ points. (The same number of
iterations is sufficient for direct solution of
the equations \tht{2.7}.) However, the above algorithm
has one clear advantage:

The function $\o(\nu\nea\l;q,t)$, which we have to compute
on each step, is much more simple than the coefficients
${\tsize \bi{\l}{\mu}{1/q,1/t}}$. In other words, one
can expand a polynomial $f$ in interpolation Macdonald 
polynomials without actually knowing those
quite complicated polynomials explicitly.

In the Jack limit this algorithm can be made
exponentially faster. Using the fact that the
LHS of \tht{5.2} depends in the Jack limit only
on $|\mu|$, not on particular form of $\mu$,
one is able to use on each step some 
previously computed sums   as
follows.

We shall use the parameter $\th$
as in \cite{OO2}. Let
$$
f(x_1,x_2,\dots)
$$
be a polynomial, symmetric in variables $x_i-\th i$.
By definition, set 
$$
\bi{\l}{\mu}{\th}=\lim_{q\to 1} \bi{\l}{\mu}{q,q^\th} \,.
$$
The polynomials
$$
\bi{x}{\mu}{\th}
$$
form a linear basis in the algebra of polynomials 
symmetric in variables $x_i-\th i$. The coefficients
$\fh(\mu)$ in the expansion
$$
f=\sum_\mu \fh(\mu) \, \bi{x}{\mu}{\th}
$$
can be found from the following specialization
of the formula \tht{2.8}
$$
\fh(\l)=\sum_{\mu\subset\l} (-1)^{|\l/\mu|} 
\bi{\l}{\mu}{\th} \, f(\mu) \,.
$$

Now, instead of computing the numbers
$$
F(\mu,\l)=(-1)^{|\l/\mu|}
\bi{\l}{\mu}{\th}\, f(\mu) \,,
$$
we shall deal with the sums
$$
\align
F^{(k)}(\l)&=\sum\Sb
\mu\subset\l\\
|\l/\mu|\le k\endSb
\binom{|\l|}{|\mu|}^{-1}
\binom{k}{|\l/\mu|} F(\mu,\l)\,,
\\
&= 
\sum\Sb
\mu\subset\l\\
|\l/\mu|\le k\endSb
(-1)^{|\l/\mu|} 
\binom{|\l|}{|\mu|}^{-1}
\binom{k}{|\l/\mu|}
\bi{\l}{\mu}{\th} \, f(\mu)\,,  
\endalign 
$$
where
$$
k=0,\dots,|\l|\,.
$$
Then we have  
$$
\align
F^{(0)}(\l)&=f(\l)\,,\\
F^{(|\l|)}(\l)&=\fh(\l)\,,\\
\endalign
$$
and the following reccurence relation 
$$
F^{(k+1)}(\l)=F^{(k)}(\l)-\frac1{|\l|}\sum_{\nu\nea\l}
\o(\nu\nea\l;\th) \, F^{(k)}(\nu)\,,
$$
where (see Fig.~2)
$$
\align
\o(\nu\nea\l;\th)&=\lim_{q\to 1}\o(\nu\nea\l;q,q^\th) \\
&=\frac{
\dsize \prod_{\text{inner corners $\bullet$ \,}}
\left(\Delta y(\bigstar,\bullet)+ \th \Delta x(\bullet, \bigstar)\right)}
{
\dsize \prod_{\text{outer corners $\circ$}}
\left(\Delta y(\bigstar,\circ) + \th \Delta x(\circ, \bigstar)\right)}\,.
\endalign
$$
The advantage of this algorithm is that 
instead of a cycle over all pairs of partitions
$\mu\subset\l$ one has a
cycle over all partition $\l$ and all {\it numbers} $k$,
such that $k\le|\l|$, which is much faster. 

In the case of one variable this algorithm specializes
to the classical Newton's algorithm for interpolation
with knots $0,1,2,\dots$ because
$$
\o((l-1)\nea(l);\th)=l \,.
$$


\Refs

\widestnumber\key{KOO}

\ref
\key Ab
\by A.~Abderrezzak
\paper G\'en\'eralisation de la transformation d'Euler
d'une s\'erie formelle
\jour Adv.\ Math.\
\vol 103 \yr 1994 \issue 2 \pages 180--195
\endref

\ref
\key An
\by G.~Andrews
\book The theory of partitions. Encyclopedia of 
Mathematics and its Applications, Vol. 2.\
\publ Addison-Wesley Publishing Co.\
\publaddr Reading, Mass.-London-Amsterdam
\yr 1976
\endref

\ref
\key Bi
\by C.~Bingham
\paper  An identity involving partitional 
generalized  binomial coefficients
\jour  J. Multivariate Analysis 
\vol 4 \yr 1974 \pages 210--223
\endref

\ref
\key Ch
\by I.~Cherednik
\paper Macdonald's evaluation conjectures and
difference Fourier transform 
\jour Invent.\ Math.\ 
\vol 122 \yr 1995 \issue 1 \pages 119--145
\endref


\ref
\key EK
\by P.~Etingof and A.~A.~Kirillov, Jr.\
\paper Representation-theoretic proof of the inner product
and symmetry identities for Macdonald's polynomials
\jour Compositio Math.\
\vol 102 \pages 179--202 \yr 1996
\endref

\ref
\key FK 
\by J.~Faraut and A. ~Kor\'anyi
\book Analysis on symmetric cones
\publ Oxford Univ.\ Press
\yr 1994
\endref

\ref
\key GR
\by G.~Gasper and M.~Rahman
\book Basic hypergeometric series
\publ Cambridge University Press
\yr 1990
\endref



\ref
\key Ka
\by J.~Kaneko
\paper Selberg integrals and hypergeometric
functions associated with Jack polynomials
\jour SIAM J.\ Math.\ Anal.\ 
\vol 24 \yr 1993 \pages 1086--1110.
\endref

\ref
\key Ka2
\bysame
\paper $q$-Selberg integrals and Macdonald  polynomials
\jour Ann.\ Sci.\ \'Ecole Norm.\ Sup.\ (4)
\vol 29 \yr 1996 \issue 5 \pages 583--637
\endref



\ref
\key Ki
\by A.~A.~Kirillov, Jr.\
\paper On inner product in modular tensor categories I
\paperinfo to appear in Journal of AMS, q-alg/9508017
\endref

\ref
\key KOO
\by S.~Kerov, A.~Okounkov, and G.~Olshanski
\paper The boundary of Young graph with Jack edge multiplicities
\paperinfo to appear in IMRN, q-alg/9703037
\endref



\ref
\key Kn
\by F.~Knop
\paper Symmetric and non-Symmetric quantum Capelli
polynomials
\paperinfo to appear
\endref


\ref
\key KS
\by F.~Knop and S.~Sahi
\paper Difference equations and symmetric polynomials
defined by their zeros
\jour Internat.\ Math.\ Res.\ Notices 
\yr 1996 \issue 10 \pages 473--486
\endref 



\ref
\key Lasc
\by A.~Lascoux
\paper Classes de Chern d'un produit tensoriel
\jour Comptes Rendus Acad.\ Sci.\ Paris, S\'er.\ I
\vol 286A \yr 1978 \pages 385--387
\endref

\ref
\key La
\by M.~Lassalle
\paper Une formule de bin\^ome
g\'en\'eralis\'ee pour les polyn\^omes de Jack
\jour Comptes Rendus
Acad.\ Sci.\ Paris, S\'er.\ I
\vol 310 \yr 1990
\pages 253--256
\endref


\ref
\key M
\by I.~G.~Macdonald
\book Symmetric functions and Hall polynomials, 
second edition
\publ Oxford University Press \yr 1995
\endref

\ref
\key Ok1
\by A.~Okounkov
\paper
Quantum immanants and higher Capelli identities
\jour Transformation groups
\vol 1 \issue 1 \yr 1996 \pages 99-126
\endref


\ref
\key Ok2
\bysame
\paper
Young Basis, Wick Formula, and Higher Capelli
identities
\jour Internat.\ Math.\ Res.\ Not.\
\vol 17 \yr 1996 \pages 817--839 
\endref


\ref
\key Ok3
\bysame
\paper
(Shifted) Macdonald polynomials: $q$-Integral
Representation and Combinatorial formula
\paperinfo
to appear in Comp.\ Math., q-alg/9605013
\endref

\ref
\key Ok4
\bysame
\paper
$BC_n$-type shifted Macdonald polynomials and binomial formula for
Koornwinder polynomials 
\paperinfo to appear, q-alg/9611011
\endref

\ref
\key Ok5
\bysame
\paper On n-point correlations in the log-gas at rational temperature
\paperinfo to appear, hep-th/9702001
\endref

\ref
\key Ok6
\bysame
\paper Proof of a conjecture of Goulden and Jackson 
\paperinfo to appear in Can.\ J.\ of Math.\ 
\endref

\ref
\key Ok7
\bysame
\paper A characterization of interpolation Macdonald polynomials
\paperinfo preprint (1997) 
\endref 

\ref
\key OO
\by A.~Okounkov and G.~Olshanski
\paper Shifted Schur functions
\jour Algebra i Analiz
\vol 9
\yr 1997
\pages No.~2
\lang Russian
\transl\nofrills English version to appear in St.~Petersburg Math. J. 
{\bf 9} (1998), No.~2
\endref

\ref
\key OO2
\bysame
\paper Shifted Jack polynomials, binomial formula,
and applications
\jour Math.\ Res.\ Letters
\vol 4 \yr 1997 \pages 69--78
\paperinfo q-alg/9608020
\endref

\ref
\key OO3
\bysame
\paper Shifted Schur functions II
\paperinfo to appear in  A.~A.~Kirillov
representation theory seminar, G.~Olshanski ed.,
Adv.\ in Math.\ Sciences, Amer.\ Math.\ Soc.
\endref

\ref 
\key OO4
\bysame
\paper Asymptotics of Jack polynomials as $\dim\to\infty$
\paperinfo in preparation
\endref

\ref
\key S1
\by S.~Sahi
\paper The Spectrum of Certain Invariant Differential Operators
Associated to a Hermitian Symmetric Space
\inbook Lie theory and geometry: in honor of Bertram Kostant,
Progress in Mathematics
\vol 123
\eds J.-L.~Brylinski, R. Brylinski, V.~Guillemin, V. Kac
\publ Birkh\"auser
\publaddr Boston, Basel
\yr 1994
\endref


\ref
\key S2
\bysame
\paper Interpolation, Integrality, and a generalization
of Macdonald's polynomials
\jour Internat.\ Math.\ Res.\ Notices 
\yr 1996 \issue 10 \pages 457--471
\endref



\endRefs
\enddocument

\end